\documentclass[12pt]{article}
\usepackage{latexsym}
\usepackage{amsfonts}
\usepackage{amsmath}
\usepackage{amsthm}
\usepackage{geometry} 
\usepackage{xcolor} 
\usepackage{amssymb}
\usepackage{parskip}
\usepackage{graphicx}
\usepackage{enumerate}
\usepackage{hyperref}
\usepackage{xspace,colortbl}

\theoremstyle{plain}
 \newtheorem{lemma}{{\bf Lemma}}

\newtheorem{remark}{{\bf Remark}}
  
 \newtheorem{theorem}[lemma]{{\sc Theorem}}

\newcommand{\be}{\begin{equation}}
\newcommand{\ee}{\end{equation}}
\newcommand{\ben}{\begin{eqnarray*}}
\newcommand{\een}{\end{eqnarray*}}

\newcommand{\bt}{\begin{theorem}}
\newcommand{\et}{\end{theorem}}
\newcommand{\bl}{\begin{lemma}}
\newcommand{\el}{\end{lemma}}
\newcommand{\bea}{\begin{eqnarray}}
\newcommand{\eea}{\end{eqnarray}}
 \usepackage{xspace,colortbl}
\definecolor{blu}{rgb}{.01,.01,1}
\definecolor{dblu}{rgb}{.01,.01,.85}
\definecolor{dred}{rgb}{.01,.01,.01}
\definecolor{red}{rgb}{1,1,1}
\definecolor{grn}{rgb}{.01,.5,.01}
 
 \newcommand{\cdr}{\color{dred}}

\newcommand{\cla}{{\mathcal A}}

\newcommand{\clf}{{\mathcal F}}

\newcommand{\xa}{{\mbox{\large\color{dred} ${\bf a}$}}}
\newcommand{\xb}{{\mbox{\large\color{dred} ${\bf b}$}}}

\newcommand{\xu}{{\mbox{\large\color{dred} ${\bf u}$}}}
\newcommand{\xv}{{\mbox{\large\color{dred} ${\bf v}$}}}

\newcommand{\CS}{{\mbox{\color{dred}${\mathbf S}$}}}

\newcommand{\tx}{{\mbox{\large\color{dred} ${\bf x}$}}}
\newcommand{\ty}{{\mbox{\large\color{dred} ${\bf y}$}}}

\newcommand{\B}{{\mathbb B}}

\newcommand{\E}{{\mathbb E}}

\renewcommand{\L}{{\mathcal L}}

\renewcommand{\P}{{\mathbb P}}

\newcommand{\R}{{\mathbb R}}

\newcommand{\Ind}{{\mathbf 1}}

\newcommand{\Bb}{{\large {\bf\mathbb B}}}

\newcommand{\Cb}{{\large {\bf\mathbb C}}}
\newcommand{\Db}{{\large {\bf\mathbb D}}}
\newcommand{\Fx}{{\mathbf F}}
\newcommand{\Gx}{{\mathbf G}}

\newcommand{\CW}{{\color{dred}{\mathbf W}}}
\newcommand{\CX}{{\color{dred}{\mathbf X}}}
\newcommand{\CY}{{\color{dred}{\mathbf Y}}}
\newcommand{\CZ}{{\color{dred}{\mathbf Z}}}

\newcommand{\Lb}{{\color{dred}{\bf L}}}
\newcommand{\xx}{{\large\color{dred}{\bf x}}}
\newcommand{\xy}{{\large\color{dred}{\bf y}}}

\newcommand{\xw}{{\large\color{dred}{\bf w}}}
\newcommand{\xz}{{\large\color{dred}{\bf z}}}

\newcommand{\nrm}[1]{\Vert\, #1\,\Vert}
\newcommand{\nrmB}[1]{\Vert\, #1\,\Vert_{_B}}

\usepackage{geometry} 
\begin{document} 
\thispagestyle{empty}
\vspace*{15mm}
\centerline{\bf \Large Stochastic approximation in infinite dimensions}
\vskip 12mm

\centerline{ \large Rajeeva L Karandikar}
\vskip 1mm
\centerline{and}
\vskip 1mm
\centerline{ \large B. V. Rao}
\vskip 6mm

\centerline{\large  In honour of K R Parthasarathy}
\vskip 11mm

\centerline{\em Chennai Mathematical Institute}
\centerline{\em H1 Sipcot IT Park, Siruseri, TN 603103, India.} 

\begin{abstract}
Stochastic Approximation (SA) was introduced in the early 1950's and has been an active area of research for several decades.  While the initial focus was on statistical questions, it was seen to have applications to signal processing, convex optimisation. 
In later years SA has found application in Reinforced Learning (RL) and led to revival of interest. 

While bulk of the literature is on SA for the case when the observations are from a  finite dimensional Euclidian space, there has been interest in extending the same to infinite dimension. Extension to Hilbert spaces is relatively easier to do, but this is not so when we come to a Banach space - since in the case of a Banach space, even {\em law of large numbers} is not true in general. We consider some  cases where approximation works in a Banach space. Our framework includes case when the Banach space $\Bb$  is $\Cb([0,1],\R^d)$, as well as $\L^1([0,1],\R^d)$, the two cases which do not even have the Radon-Nikodym property. 
\end{abstract}

\vfill

\hrule

\vspace*{3mm}


\noindent
{\em email: rlk@cmi.ac.in, \hskip 2mm bvrao@cmi.ac.in}
\newpage

\section{Introduction}

Let $f:\R^d\mapsto \R^d$ be a function such that the equation $f(x)=0$ has a unique solution, say $x^*$.  If  one could observe $f(x)$ given $x$, then the solution $x^*$ can be approximated by a suitable successive approximation procedure, the exact algorithm would depend upon the conditions on $f$.  Now if $f(x)$ is not directly observable but can be observed only in presence of noise, then  Robbins and Monro  introduced a method, now known as Stochastic Approximation (SA), to approximate the solution $x^*$: Choose $x_0\in\R$ and for $n\ge 1$ define
\[
x_{n+1}= x_n+ \alpha_n(f(x_n)+\xi_{n+1}),\;\;n\ge 1,\]
where $\xi_{n+1}$ is the additive observational noise and $\alpha_n\in(0,1)$ is the {\it step size sequence} satisfying 
\begin{equation}\label{k1}\textstyle \sum_k\alpha_k=\infty \text{ and } \sum_k\alpha^2_k<\infty,\end{equation}  then $\{x_{n+1}:n\ge 0\}$ converges to the true solution $x^*$. While Robbins and Monro \cite{RM}  had considered $d=1$, $f$ to be a monotone function, $\{\xi_{n+1}:n\ge 0\}$  to be a sequence of independent random variables satisfying appropriate moment conditions and the convergence of $x_{n+1}$ to $x^*$ was asserted only in mean square, the same was extended in multiple directions and has been applied in several contexts. See the {\em Introduction} in  \cite{Borkar2008} and also, \cite{BMP92}, \cite{Bertsekas}, \cite{Blum54}, \cite{Borkar98}, \cite{Borkar-Meyn00}, \cite{Jaak}, \cite{KV1}, \cite{KV2}, \cite{LT}, \cite{Lai}, \cite{Milz}. 

While most extensions were in the case when the state space is finite dimensional, there have been some extension to the case when the state space is an Hilbert space or a Banach space. As can be expected, the analysis in $\R^d$ case can be extended to Hilbert spaces with some work (see \cite{Dieuleveut}, \cite{KS}), the case of Banach Spaces is more complicated as even the law of large numbers may not be valid on a Banach space (see \cite{FJPP}, \cite{LT})
unless we impose some conditions on geometry of the underlying Banach space. 

In this paper, we consider the case when the state space is a Banach space $\Bb$. 
Let $\Gx:\Bb\mapsto \Bb$ be a function such that $\Gx(x)=0$ has a unique solution. For example, $\Gx$ could be $\Fx(\xx)-\xx$ where $\Fx$ is a contraction or $t\mapsto \Gx(\xx+t\xy)$ is a strictly monotone function for all $\xx,\,\xy\in\Bb$. Let the  unique solution to $\Gx(\xx)=0$ be denoted by $\xx^*$. Then for any $\xu_0\in\Bb$, the sequence
\[\xu_{n+1}=\xu_n+\Gx(\xu_n),\;\;n\ge 0\]
converges to $\xx^*$ (under appropriate conditions on $\Gx$). Now suppose that $\Gx(\xu_n)$ is not directly observable, but can only be observe corrupted by (additive) noise, {\em i.e.} - we observe $\Gx(\xu_n)+\CZ_{n+1}$, where $\CZ_{n+1}$ is noise. One can instead modify the approximation algorithm as follows (taking cue from the finite dimensional case): for a sequence of constants $\alpha_n\in (0,1)$ converging to zero, we define
\[\CX_{n+1}=\CX_n+\alpha_n(\Gx(\CX_n)+\CZ_{n+1}),\;\;n\ge 0.\]
Under suitable conditions on the function $\Gx$, the noise sequence $\{\CZ_n:n\ge 0\}$ and the step size $\{\alpha_n: n\ge 0\}$, one hopes to conclude that $\CX_n$ converges to $\xx^*$ almost surely.

When $\Bb=\R^d$,  if one assumes that\\
 (i) the noise sequence is a martingale difference sequence with bounded variance, \\(ii) the step size sequence satisfies equation \eqref{k1}, known as the {\em Robbins-Monroe conditions}, \\(iii) $\Gx(\xx)=\Fx(\xx)-\xx$, where $\Fx$ is a contraction,  or $\Gx(\xx)=\nabla\Fx$, where  $\Fx$ is a smooth convex function,\\ the required conclusion follows  (with some additional conditions on $\Fx$). 
We will show that the same is valid in the infinite dimensional case under several sets of assumptions, including the case when $\Bb$ is an arbitrary separable Banach space,  when the noise sequence is an i.i.d. sequence with  Gaussian distribution. When $\Bb=\Cb([0,1],\R^d)$, we will show that the conclusion holds when the noise is a sequence of independent martingales with a suitable condition on the moments. 

For more general noise sequences, if the underlying space $\Bb$ is a Hilbert space or if the underlying Banach space satisfies some suitable geometric conditions, we will show that the result holds under fairly general conditions on the noise sequence analogous to the Kolmogorov three series theorem (see Brown \cite{Brown}).

\section{Notations and Main Results}
Let $\Bb$ be a Banach space. Let $\Gx:\Bb\mapsto \Bb$ be a function such that the equation $\Gx(\xx)=0$ for $\xx\in\Bb$ admits a unique solution  $\xx^*$. 
Here and in the rest of the paper,  $\nrmB{\xx}$ denotes the norm, for $\xx\in\Bb$. 
Let us start with some scenarios that we would like to include in our analysis.

\noindent{\bf Examples:} 
\begin{enumerate}[a)]
\item Let $\Fx:\Bb\mapsto \Bb$ be a contraction:
\begin{equation}\label{R1}
\nrmB{\Fx(\xx)-\Fx(\xy)}\le \gamma\nrmB{\xx-\xy},\;\;\forall \xx, \xy\in\Bb,\;0<\gamma<1
\end{equation}
with (unique) fixed point $\xx^*$: $\Fx(\xx^*)=\xx^*$. Let us note that \eqref{R1} implies
\begin{equation}\label{t3}
\nrmB{\Fx(\xx)-\xx^*}\le \gamma\nrmB{\xx-\xx^*}.
\end{equation}
and then writing $\Gx(\xx)=\xx-\Fx(\xx)$, we see that $\Gx(\xx)-(\xx-\xx^*)=\Fx(\xx)-\Fx(\xx^*)$ and hence condition \eqref{R1} can be rewritten as 
\begin{equation}\label{R0}
\nrmB{\Gx(\xx)-(\xx-\xx^*)}\le \gamma\nrmB{\xx-\xx^*}.
\end{equation} 
\item Let $\Bb$ be a Hilbert space and $\Fx:\Bb\mapsto \R$ be a $C^2$-convex function { and let  $\Gx=\nabla \Fx$. Suppose that the eigenvalues of Hessian $\Delta \Fx$ are bounded below by $c_1>0$ and bounded above by $c_2<\infty$. Let $\xx^*$ denote the unique solution to $\Gx(\xx)=0$. Note that the dual  $\Bb^*=\Bb$ and denoting the action of $\xu\in\Bb^*=\Bb$ on $\xx\in\Bb$ by $\langle \xu,\xx\rangle$,  for any $\xx,\xy, \xu \in \Bb$, we have
\begin{equation}\label{S3}
\langle\xu, \xx-\xy\rangle\langle\xu, \Gx(\xx)-\Gx(\xy)\rangle\ge 0
\end{equation}
\begin{equation}\label{S4}
c_1|\langle\xu,\xx-\xy\rangle|  \le|\langle\xu,\Gx(\xx)-\Gx(\xy)\rangle|\le c_2|\langle\xu,\xx-\xy\rangle|.
\end{equation} 
Let us note that \eqref{S3} and \eqref{S4} imply that with $\theta={1+c_1+c_2}$ and $\rho=1- \frac{c_1}{1+c_1+c_2}$, for any $\xu\in\B^*=\Bb$, 
$\langle\xu,\theta^{-1}\Gx(\xx)\rangle$ has same sign as $\langle\xu,\xx-\xx^*\rangle$ and  lies between $\frac{c_1}{1+c_1+c_2}\langle\xu,\xx-\xx^*\rangle$ and  $\frac{c_2}{1+c_1+c_2}\langle\xu,\xx-\xx^*\rangle$. Hence
\[|\langle\xu,(\xx-\xx^*)-\theta^{-1}\Gx(\xx)\rangle|\le \rho |\langle\xu,\xx-\xx^*\rangle|.\]
Taking supremum over $\xu\in\B^*$, we conclude that 
\begin{equation}\label{S9}\nrmB{\theta^{-1}\Gx(\xx)-(\xx-\xx^*)}\le \rho\nrmB{\xx-\xx^*}.\end{equation}
\item Let $\Gx:\Bb\mapsto \Bb$ be a function such that the equation $\Gx(\xx)=0$ for $\xx\in\Bb$ has a unique solution  $\xx^*$ and 
$\exists$ constants $0<c_1<c_2<\infty$ such that  $\forall \xx, \xy\in\Bb$, $\forall \Lb\in\Bb^*$
on has
\begin{equation}\label{R3}
(\Lb(\xx-\xy))(\Lb(\Gx(\xx)-\Gx(\xy)))\ge 0
\end{equation}
\begin{equation}\label{R4}
c_1|\Lb(\xx-\xy)|  \le|\Lb(\Gx(\xx)-\Gx(\xy))|\le c_2|\Lb(\xx-\xy)|.
\end{equation} 
Let us note that \eqref{R3} and \eqref{R4} imply that for $\theta={1+c_1+c_2}$ and $\rho=1- \frac{c_1}{1+c_1+c_2}$, for any $\Lb\in\B^*$, 
$\Lb(\theta^{-1}\Gx(\xx))$ has same sign as $\Lb(\xx-\xx^*)$ and  lies between $\frac{c_1}{1+c_1+c_2}\Lb(\xx-\xx^*)$ and  $\frac{c_2}{1+c_1+c_2}\Lb(\xx-\xx^*)$
hence
\[|\Lb((\xx-\xx^*)-\theta^{-1}\Gx(\xx)|\le \rho |\Lb(\xx-\xx^*)|.\]
Taking supremum over $\Lb\in\B^*$, we get \eqref{S9}}
\end{enumerate}

Thus we will assume that the function  $\Gx:\Bb\mapsto \Bb$  and that the equation $\Gx(\xx)=0$ for $\xx\in\Bb$ has a unique solution  $\xx^*$ and that $\exists$ constants $\rho,\theta$, $0<\rho<1$ and $1\le\theta<\infty$ such that 
\begin{equation}\label{R2}
\nrmB{\theta^{-1}\Gx(\xx)-(\xx-\xx^*)}\le \rho\nrmB{\xx-\xx^*}\;\;\;\forall \xx\in\Bb.
\end{equation}

All random variables considered in this paper will be defined on one probability space $(\Omega,\clf,\P)$.
Throughout this paper, $\{\alpha_n:n\ge 0\}$ denotes $(0,1)$-valued sequence of constants that satisfies 
\begin{equation}\label{R5}
\textstyle \lim_n\alpha_n=0,\end{equation}
\begin{equation}\label{R22}
\textstyle\sum_{n=0}^\infty\alpha_n=\infty
\end{equation}
and a function $\Gx:\Bb\mapsto\Bb$ satisfying \eqref{R2}. 

 {\color{dred} Here is the first result for the case when the errors are Gaussian.}
\begin{theorem}\label{T1} Suppose $\Bb$ is a separable Banach Space.
Suppose that the noise  $\{\CZ_n:n\ge 0\}$ is a sequence of independent identically distributed $\,\Bb$ valued Gaussian random variables. Suppose the step size $\{\alpha_n: n\ge 0\}$ satisfies, in addition to \eqref{R22}, 
\begin{equation}\label{3a}\textstyle\sum_{n=0}^\infty\alpha^2_n<\infty.\end{equation}
 Let $\{\CX_n: n\ge 0\}$ be defined by $\CX_0=\xx_0$ and for $n\ge 0$:
\begin{equation}\label{3b}
\CX_{n+1}=\CX_n-\alpha_n\Bigl(\Gx(\CX_n)+\CZ_{n+1}\Bigr),\;\;n\ge 0\end{equation}
Then $\{\CX_{n+1}:n\ge 0\}$ converges to $\xx^*$ almost surely.
\end{theorem}
{\color{dred} Going beyond Gaussian errors, we first consider the case when errors are martingales on a suitable path  space.} 
\begin{theorem}\label{T2}
Let $\Bb=\Cb([0,1],\R^d)$ or $\Bb=\Db([0,1],\R^d)$ with the {\bf sup} norm: for $\xu\in\Bb$
\[\nrmB{\xu}=\sup_{0\le t\le 1}|\xu_t|.\] 
 Let $\{\CZ_n:n\ge 0\}$ be a sequence of independent $\Bb$ valued martingales with mean 0 and bounded variance {\em i.e.}
for each $n$, $\CZ_{n,0}=0$, for $0\le s<t\le 1$
\begin{equation}\label{4e}
\E[\CZ_{n,t}|\sigma(\CZ_{n,u}:u\le s)]=\CZ_{n,s}
\end{equation}
\begin{equation}\label{4f}
\E[\CZ_{n,1}^2]\le \sigma^2\end{equation}
for a constant $\sigma^2<\infty$. We assume that the step size sequence satisfies
the Robbins-Monroe conditions, \eqref{R22} and \eqref{3a}. 
Then $\{\CX_{n+1}:n\ge 0\}$ (defined by \eqref{3b}) converges to $\xx^*$ almost surely. 
\end{theorem}
{ To go beyond martingale errors, we need to impose geometric conditions on the underlying Banach space. A Banach space $\Bb$  is said to be $p$-uniformly smooth, $1< p\le 2$,  if $\exists D: 0<D<\infty$ such that
\begin{equation}\label{4g}
 \nrmB{\tx+\ty}^p+\nrmB{\tx-\ty}^p\le 2\nrmB{\tx}^p+D\nrmB{\ty}^p,\;\;\forall \tx,\ty\in\Bb.
 \end{equation}
 
See Woyczynski, \cite{Woyczynski}. Clearly, a Hilbert space is $2$-uniformly smooth.
It can be verified that $\L^p([0,1],\R^d)$, $\L^p(\tilde{\Omega},\tilde{\clf},\tilde{\P})$ for any $\tilde{\P}$ are $p$-smooth if $1<p\le 2$ and are 2-uniformly smooth if $p>2$.

It may be noted that the Banach spaces
$\Cb([0,1],\R^d)$ and $\Db([0,1],\R^d)$
 do not satisfy these geometric conditions.}

In order to include the known results in the finite dimensional case, we consider a slightly more general case.  Let $\lambda_n:\Bb^{n+1}\mapsto\R$ be a function that satisfies, for some constant $C<\infty,\;\;\forall \xu_j\in \Bb,\;n\ge 0$
\begin{equation}\label{3b1}
|\lambda_n(\xu_0,\xu_1,\ldots,\xu_n)|\le C(1+\max_{0\le k\le n}\nrmB{\xu_k}).
\end{equation}

\begin{theorem}\label{T3}
Let $\Bb$ be a $p$-uniformly smooth Banach space for some $p\in (1,2]$. 

Suppose that $\Bb$ valued noise sequence $\{\CZ_n:n\ge 0\}$ satisfies, for some constant sequences $\{\alpha_k:k\ge 0\}$, $\{\mu_k:k\ge 0\}$ and $\{\sigma_k^2:k\ge 0\}$
\begin{equation}\label{4}
\textstyle\sum_{n=0}^\infty\P(\nrmB{\alpha_n\CZ_{n+1}}\ge 1)<\infty,
\end{equation}
and writing $\overline{\CZ}_{n+1}=\E[\CZ_{n+1}{\Ind_{_{\{\nrmB{ \alpha_n\CZ_{n+1}}<1\}}}}  |\sigma(\CZ_j:j\le n)]$, we have
\begin{equation}\label{4a}
\E\bigl[\nrmB{\overline{\CZ}_{n+1}}\bigr]\le \mu_n,\;\;n\ge 0,
\end{equation}
\begin{equation}\label{4b}
\E[\nrmB{\CZ_{n+1}\Ind_{_{\{\nrmB{\alpha_{n}\CZ_{n+1}}<1\}}}-\overline{\CZ}_{n+1}}^2]\le \sigma^2_n,\;\;n\ge 0
\end{equation}
where  the step size $\{\alpha_n: n\ge 0\}$ satisfies (in addition to \eqref{R5} and \eqref{R22})
\begin{equation}\label{4c1}
\textstyle\sum_{n=0}^\infty\alpha_n\mu_n<\infty,
\end{equation}
\vspace{-12mm}

and 
\vspace{-12mm}

\begin{equation}\label{4d1}
\textstyle\sum_{n=0}^\infty\alpha^p_n\sigma^p_n<\infty.
\end{equation}
Let $\CY_0=\xy_0$ and for $n\ge 0$
\begin{equation}\label{3b2}
\CY_{n+1}=\CY_n-\alpha_n\Bigl(\Gx(\CY_n)+\lambda_n(\CY_0,\ldots,\CY_n)\CZ_{n+1}\Bigr),\;\;n\ge 0\end{equation}

Then $\{\CY_{n+1}:n\ge 0\}$  converges to $\xx^*$ almost surely. 
\end{theorem}
\begin{remark}
In the finite dimension case, often in the equation \eqref{3b2}, $\lambda_n$ is taken as a function only of $\CY_n$. 
\end{remark} 
Note that a Hilbert space is also a 2-uniformly smooth Banach space. Also, it is known that 
$\L^p([0,1],\R^d)$ for $p>2$ is a 2-uniformly smooth Banach space, while  $\L^p([0,1],\R^d)$
for $1<p<2$ is a $p$-uniformly smooth Banach space.

Here is another result which covers the case $\Bb=\L^1([0,1],\R^d)$, which is not covered by Theorem \ref{T3}.
\begin{theorem}\label{T4}
Let $\Bb=\L^p([0,1],\R^d)$, $1\le p\le 2$ with $\nrmB{\xx}=\bigl(\int_0^t|\xx_t|^pdt\bigr)^\frac{1}{p}$. 

Suppose that the noise sequence $\{\CZ_n:n\ge 0\}$ satisfies, for some constant sequences $\{\delta_k:k\ge 0\}$, $\{\mu_k:k\ge 0\}$ and $\{\sigma_k^2:k\ge 0\}$
\begin{equation}\label{4p}
\E\bigl[\Bigr(\int_0^1\bigl[|\CZ_{n+1,t}|^p\Ind_{_{\{|\alpha_n\CZ_{n+1,t}|\ge 1\}}}\bigr]dt\Bigl)^\frac{1}{p}\bigr]\le\delta_n,
\end{equation}
and writing $\overline{\CZ}_{n+1,t}=\E\bigl[\CZ_{n+1,t}\Ind_{_{\{\alpha_n|\CZ_{n+1,t}|<1\}}}|\sigma(\CZ_j:j\le n)\bigr]$, we have
\begin{equation}\label{4q}
\E\bigl[\Bigr(\int_0^1|\overline{\CZ}_{n+1,t}|^pdt\Bigl)^\frac{1}{p}\bigr]\le\mu_n,
\end{equation}
\begin{equation}\label{4r}
\E\bigl[\int_0^1\bigl|\CZ_{n+1,t}\Ind_{_{\{|\alpha_n\CZ_{n+1,t}|<1\}}}-\overline{\CZ}_{n+1,t}\bigr|^2dt\bigr]\le \sigma^2_n,\;\;n\ge 0.
\end{equation}
Suppose that the step size $\{\alpha_n: n\ge 0\}$ satisfies (in addition to \eqref{R22})
\begin{align}\label{4s}
\textstyle\sum_{n=0}^\infty\alpha_n\delta_n&<\infty,\\
\label{4t}
\textstyle\sum_{n=0}^\infty\alpha_n\mu_n&<\infty,\\
\label{4u}
\textstyle\sum_{n=0}^\infty\alpha^2_n\sigma^2_n&<\infty.
\end{align}

Then $\{\CY_{n+1}:n\ge 0\}$ (defined by \eqref{3b2}) converges to $\xx^*$ almost surely. 
\end{theorem}
\begin{remark} Let us note that  $\alpha_n=\frac{1}{n\log(n)}$, $\delta_n=\mu_n=\frac{1}{\log(n)}$ and $\sigma_n^2=n$ satisfy the conditions \eqref{4s},  \eqref{4t} and  \eqref{4u}.
\end{remark}
\section{The Deterministic case}
Here is the key result that allows us to separate the probabilistic part (which is convergence of the weighted sum of the noise terms) and deterministic part- which shows that if the series of  weighted sum of the noise terms converges then approximating sequence of observations ($\{\CX_n: n\ge 0\} $ in Theorems \ref{T1} and \ref{T2} and $\{\CY_n: n\ge 0\} $ in Theorems \ref{T3} and \ref{T4}) also converge to the desired limit.

\begin{theorem}\label{T0}\label{deterministic}
Let $\{\xz_{n+1}:n\ge 0\}\subseteq \Bb$ be such that
\begin{equation}\label{R9}
\textstyle\sum_{n=0}^m\alpha_n\xz_{n+1}\;\;\;\text{converges in $\Bb$ \mbox{(in $\nrmB{\cdot}$)} as $m\rightarrow\infty$.}
\end{equation}
For $n\ge 0$, let $\psi_n:\Bb^{n+1}\mapsto\R: n\ge 0\}$ be given by
\begin{equation}\label{R12}
\psi_n(\xu_0,\xu_1,\ldots,\xu_n)=(1+\max_{0\le k\le n}\nrmB{\xu_k}).
\end{equation}
For $\xx_0\in\Bb$ fixed, and for $n\ge 0$ let 
\begin{equation}\label{R6}
\xx_{n+1}=\xx_n-\alpha_n\Gx(\xx_n)+\psi_n(\xx_0,\ldots,\xx_n)\alpha_n\xz_{n+1},\;\;n\ge 0.\end{equation}
Then $\xx_n$ converges to $\xx^*$.
\end{theorem}
{\bf\cdr Proof:} In view of \eqref{R5}, by ignoring first few terms and 
re-labbeling, we can assume that  $\beta_n<1$ for $n\ge 0$. Now \eqref{R6} can be rewritten as
\[
\begin{split}
\xx_{n+1}=&(1-\beta_n)\xx_n-\beta_n(\textstyle\theta^{-1}\Gx(\xx_n)-\xx_n)-\psi_n(\xx_0,\ldots,\xx_n)\alpha_n\xz_{n+1}\\
=&(1-\beta_n)\xx_n-\beta_n(\xy_n-\xx_n)-\phi_n\alpha_n\xz_{n+1}
\end{split}\]
where $\beta_n=\theta\alpha_n$, $\xy_n=\theta^{-1}\Gx(\xx_n)$ and $\phi_n=\psi_n(\xx_0,\ldots,\xx_n)$ and thus
\begin{equation}\label{R7}
\xx_{n+1}-\xx^*=(1-\beta_n)(\xx_n-\xx^*)-\beta_n(\xy_n-(\xx_n-\xx^*))-\phi_n\alpha_n\xz_{n+1}.
\end{equation}
Note that in view of the assumption \eqref{R2} we have
\begin{equation}\label{b25bx} 
\nrmB{(\xy_k-(\xx_k-\xx^*))}\le \rho\nrmB{\xx_k-\xx^*}\\
\end{equation}
We can check, using recursion and \eqref{R7} 
that for $ m\le n$,  we have
 \begin{equation}\label{R8}
 \begin{split}
\xx_{n+1}-\xx^*=&\textstyle\bigl[\prod_{j=m}^{n}(1-\beta_j)\bigr](\xx_m-\xx^*)\\&\;\;\;\;-\textstyle\sum_{k=m}^{n}\textstyle\bigl[\prod_{_{\{j:k<j\le n\}}}(1-\beta_j)\bigr]\beta_k(\xy_k-(\xx_k-\xx^*))\\
&\;\;\;\;-\textstyle\sum_{k=m}^{n}\textstyle\bigl[\prod_{_{\{j:k<j\le n\}}}(1-\beta_j)\bigr]\phi_k\alpha_k\xz_{k+1}. 
\end{split}\end{equation}
Here, product over an empty set is taken as 1.
 Let us introduce notation for $ m\le n$,

\vspace*{-9.8mm}
 \begin{align*}
 \xa_{m,n}&=\textstyle{\sum_{k=m}^n}\alpha_k\xz_{k+1}\\
\xb_{m,n}&=\textstyle{\sum_{k=m}^n}\bigl[\prod_{\{k<j\le n\}}(1-\beta_j)\bigr]\alpha_k\xz_{k+1}\\
 \xu_{m,n}&=\textstyle\bigl[\prod_{j=m}^{n}(1-\beta_j)\bigr](\xx_m-\xx^*)\\
 \xv_{m,n}&=\textstyle\sum_{k=m}^{n}\textstyle\bigl[\prod_{_{\{j:k<j\le n\}}}(1-\beta_j)\bigr]\beta_k(\xy_k-(\xx_k-\xx^*))\\
\xw_{m,n}&=\textstyle{\sum_{k=m}^n}\bigl[\prod_{\{k<j\le n\}}(1-\beta_j)\bigr]\phi_k\alpha_k\xz_{k+1}
\end{align*}
\vspace{-6mm}

In view of \eqref{R8}, we have for any $m\le n$
\[\xx_{n+1}-\xx^*=\xu_{m,n}-\xv_{m,n}-\xw_{m,n}\]
and hence
 \begin{equation}\label{R19}
\nrmB{\xx_{n+1}-\xx^*}\le \nrmB{\xu_{m,n}}+\nrmB{\xv_{m,n}}+\nrmB{\xw_{m,n}}.
\end{equation}
Let $\gamma_n=\max\{\nrmB{\xx_k-\xx^*},\;\;0\le k\le n\}$ and $\gamma^*=\sup_n \gamma_n$. Then $\nrmB{\xx_n}\le (\gamma_n+\nrmB{\xx^*})$ and hence in view of assumption \eqref{R12}, it follows that
\begin{equation}\label{R25}
\phi_n\le (1+\gamma_n+\nrmB{\xx^*}),\;\;\forall n.\end{equation}
We will first prove that for all
$\epsilon>0$, $\exists$ $n_\epsilon $ such that 
\begin{equation}\label{R24}
\nrmB{\xw_{m,n}}\le(1+\gamma_n+\nrmB{\xx^*})\epsilon\;\;\;\forall n\ge m\ge n_\epsilon.
\end{equation}
 In view of assumption \eqref{R9}, it follows that $\forall \epsilon>0$, $\exists$ $n_\epsilon$ such that 
\begin{equation}\label{R10}
\nrmB{\xa_{m,n}}\le \textstyle\epsilon\;\;\;\forall n\ge m\ge n_\epsilon.
\end{equation} 
Note that for all  $ m< n$, 
\[
\begin{split}
\xb_{m,n}&=\textstyle{\sum_{k=m}^n}[\prod_{\{k<j\le n\}}(1-\beta_j)]\alpha_k\xz_{k+1}\\
&=\textstyle{\sum_{k=m}^n}[\prod_{\{k<j\le n\}}(1-\beta_j)]({\xa}_{k,n}-{\xa}_{k+1,n})\\
\end{split}\]
where ${\xa}_{n+1,n}=0$. As a consequence, 
\[\begin{split}
\xb_{m,n}=&\textstyle{\sum_{t=m+1}^n}[\prod_{\{t<j\le n\}}(1-\beta_j)-\prod_{\{t-1<j\le n\}}(1-\beta_j)]
\xa_{t,n}+\\
&+ \textstyle[\prod_{\{m<j\le n\}}(1-\beta_j)]\xa_{m,n}\\
\end{split}\]
Thus, in view of \eqref{R10}, we have $\forall n\ge m\ge n_\epsilon$,
\be\label{R11}
\begin{split}
\nrmB{\xb_{m,n}}\le& \textstyle{\sum_{t=m+1}^n}[\prod_{\{t<j\le n\}}(1-\beta_j)-\prod_{\{t-1<j\le n\}}(1-\beta_j)]
\nrmB{\xa_{t,n}}\\
&\vspace{9mm}+ \textstyle[\prod_{\{m<j\le n\}}(1-\beta_j)]\nrmB{\xa_{m,n}}\\
\le& \textstyle{\sum_{t=m+1}^n}[\prod_{\{t<j\le n\}}(1-\beta_j)-\prod_{\{t-1<j\le n\}}(1-\beta_j)]
 \textstyle\epsilon+\\
&+ \textstyle[\prod_{\{m<j\le n\}}(1-\beta_j)] \textstyle\epsilon\\
\le &  \textstyle\epsilon.
\end{split}\ee
Now, noting that
\[
\begin{split}
\xw_{m,n}=&\textstyle{\sum_{k=m}^n}\phi_k[\prod_{\{k<j\le n\}}(1-\beta_j)]\alpha_k\xz_{k+1}\\
=&\phi_m\textstyle{\sum_{k=m}^n}[\prod_{\{k<j\le n\}}(1-\beta_j)]\alpha_k\xz_{k+1}\\&+\textstyle{\sum_{k=m+1}^n}(\phi_k-\phi_m)[\prod_{\{k<j\le n\}}(1-\beta_j)]\alpha_k\xz_{k+1}\\
=&\phi_m\xb_{m,n}+\textstyle{\sum_{k=m+1}^n}(\phi_k-\phi_m)[\prod_{\{k<j\le n\}}(1-\beta_j)]\alpha_k\xz_{k+1}\\
\end{split}\]
Proceeding likewise, we will get
\[\xw_{m,n}=\phi_m\xb_{m,n}+(\phi_{m+1}-\phi_m)\xb_{m+1,n}+\ldots +(\phi_{n}-\phi_{n-1})\xb_{n,n}\]
Since $\phi_k$ is increasing, in view of \eqref{R11} it follows that $\forall n\ge m\ge n_\epsilon$,
\[\begin{split} \nrmB{\xw_{m,n}}&\le \phi_m\nrmB{\xb_{m,n}}+(\phi_{m+1}-\phi_m)\nrmB{\xb_{m+1,n}}+\ldots +(\phi_{n}-\phi_{n-1})\nrmB{\xb_{n,n}}\\
&\le \phi_m \textstyle\epsilon+(\phi_{m+1}-\phi_m) \textstyle\epsilon+\ldots +(\phi_{n}-\phi_{n-1}) \textstyle\epsilon\\
&=\phi_{n} \textstyle\epsilon\\
&\le(1+\gamma_n+\nrmB{\xx^*})\epsilon\end{split}
\]
where we have used \eqref{R25} at the last step. 
This proves \eqref{R24}. Next we will prove that
\begin{equation}\label{b21}
\gamma^*=\sup\{\nrmB{\xx_n-\xx^*}:n\ge 0\}<\infty.
\end{equation}
 First, using $\rho<1$, get $\epsilon>0$   such that  $\epsilon<\frac{1}{2}$
\begin{equation}\label{R14}
\rho(1+2\epsilon)\le 1.
\end{equation}
In what follows, $\epsilon$ and $n_\epsilon$ are fixed so that \eqref{R14} and \eqref{R24} hold.   We will prove, by induction that $\forall k\ge n_\epsilon$
\begin{equation}\label{R13}
\gamma_k\le (1+2\epsilon)(1+\gamma_n+\nrmB{\xx^*}).
\end{equation}
Clearly \eqref{R13} is true for $k=n_\epsilon$.   
Suppose that \eqref{R13} holds for all $k$ such that  $n_\epsilon\le k\le n^*$.
We will prove that \eqref{R13} also holds for $k=n^*+1$, completing the proof that \eqref{R13} is true for all $k\ge n_\epsilon$.  To prove the induction step, suffices to prove  that
\[\nrmB{\xx_{n^*+1}-\xx^*}\le (1+2\epsilon)(1+\gamma_n+\nrmB{\xx^*}).\]

For $n_\epsilon\le k\le n^*$, \eqref{b25bx}, \eqref{R14}, \eqref{R13} imply that
\begin{equation}\label{b25a} \begin{split}
\nrmB{(\xy_k-(\xx_k-\xx^*))}&\le \rho\nrmB{\xx_k-\xx^*}\\
&\le \rho\gamma_k\\
&\le \rho (1+2\epsilon)(1+\gamma_n+\nrmB{\xx^*})\\
&\le (1+\gamma_n+\nrmB{\xx^*}).
\end{split}
\end{equation}
As a consequence
\begin{equation}\label{R15}\begin{split}
\nrmB{\xv_{n_\epsilon,n^*}}&\le \textstyle\sum_{k=n_\epsilon}^{n^*}\textstyle\bigl[\prod_{_{\{j:k<j\le n\}}}(1-\beta_j)\bigr]\beta_k\nrmB{(\textstyle\xy_k-(\xx_k-\xx^*))}\\
&\le \gamma_{n_\epsilon}\textstyle\sum_{k=n_\epsilon}^{n^*}\textstyle\bigl[\prod_{_{\{j:k<j\le n^*\}}}(1-\beta_j)\bigr]\beta_k.
\end{split}\end{equation}
It is easy to see that
\begin{equation}\label{b25b}\begin{split}
 \nrmB{\xu_{n_\epsilon,n^*}}&= \textstyle\bigl[\prod_{j=n_\epsilon}^{n^*}(1-\beta_j)\bigr]\nrmB{\xx_{n_\epsilon}-\xx^*}\\
   &\le \gamma_{n_\epsilon}\textstyle\bigl[\prod_{j=n_\epsilon}^{n^*}(1-\beta_j)\bigr].\end{split}
 \end{equation}
For all $m\le n$ one can verify (by induction on $n$) that
\begin{equation}\label{R17}
\textstyle\bigl[\prod_{k=m}^{n}(1-\beta_k)\bigr]+\textstyle\sum_{k=m}^{n}\textstyle\bigl[\prod_{_{\{j:k<j\le n\}}}(1-\beta_j)\bigr]\beta_k=1,    
\end{equation}
and then putting together \eqref{R15}, \eqref{b25b}, \eqref{R17} we conclude that
\begin{equation}\label{R20}
\nrmB{\xu_{n_\epsilon,n}}+\nrmB{\xv_{n_\epsilon,n}}\le (1+\gamma_n+\nrmB{\xx^*}).\end{equation}
Using \eqref{R24}, we deduce
   \begin{equation}\label{R21}
\begin{split}
\nrmB{\xw_{n_\epsilon,n^*}}
&\le \epsilon(\gamma_{n^*}+1+\nrmB{\xx^*})\\
&\le \epsilon(1+2\epsilon)(1+\gamma_n+\nrmB{\xx^*})\\
&\le 2\epsilon(1+\gamma_n+\nrmB{\xx^*}). 
\end{split}
\end{equation}
We thus conclude, using \eqref{R19}, \eqref{R20} along with \eqref{R21}, that
\begin{equation}\label{b31}\begin{split}
\nrmB{\xx_{n^*+1}-\xx^*}&\le \nrmB{\xu_{n_\epsilon,n^*}}+\nrmB{\xv_{n_\epsilon,n^*}}+\nrmB{\xw_{n_\epsilon,n^*}}\\
&\le (1+2\epsilon)(1+\gamma_n+\nrmB{\xx^*}).
\end{split}\end{equation}
This completes the induction proof and thus \eqref{R13} holds for all $k\ge n_{\epsilon}$  proving \eqref{b21}.

It also follows that $\sup\{\nrmB{\xx_k}: k\ge 0\}<\infty $ and hence that 
\begin{equation}\label{b51}
\sup\{\phi_{n}:n\ge 0\}<\infty.
\end{equation}
Using \eqref{R24} and \eqref{b51}, it follows that for all $\delta>0$, we can get $m_\delta$ such that
 \begin{equation}\label{b34}  
  \nrmB{\xw_{m,n}}<\delta \;\;\;\forall m,n \text{ such that }m_\delta\le m\le n
  \end{equation}
and thus we have for $m\ge m_\delta$
\begin{equation}\label{b35}
\limsup_{n\rightarrow\infty}\nrmB{\xw_{m,n}}\le\delta.
\end{equation}

Note that in view of the assumption \eqref{R22} and $\beta_n=\theta\alpha_n$, we have
\begin{equation}\label{b32}
   \lim_{n\rightarrow\infty}\bigl[\prod_{j=m}^{n}(1-\beta_j)\bigr]=0, \;\;\;\forall m<\infty.
  \end{equation} 
For $m\ge 1$, let $\zeta_m=\sup\{\nrmB{\xx_k-\xx^*}: k\ge m\}$. Then clearly, $\zeta_m$ is decreasing and $\zeta_m\le \gamma^*<\infty$ for all $m$. Let $\zeta^*=\lim_{m\rightarrow\infty}\zeta_m$. Clearly $\zeta^*\le \gamma^*$. Remains to show that $\zeta^*=0$.
Since
\[ \nrmB{\xu_{m,n}}\le \    =\nrmB{\xx_m-\xx^*}\textstyle\bigl[\prod_{j=m}^{n}(1-\beta_j)\bigr]\]
and thus in view of \eqref{b32} we have
\begin{equation}\label{b33}
\lim_{n\rightarrow\infty}\nrmB{\xu_{m,n}}=0 \;\;\;\forall m<\infty.
  \end{equation} 
Note that for any $m\le n$ we have
\begin{equation}\label{b37}\begin{split}
\nrmB{\xv_{m,n}}&\le \textstyle\sum_{k=m}^{n}\textstyle\bigl[\prod_{_{\{j:k<j\le n\}}}(1-\beta_j)\bigr]\beta_k\nrmB{(\textstyle\xy_k-(\xx_k-\xx^*))}\\
&\le \textstyle\sum_{k=m}^{n}\textstyle\bigl[\prod_{_{\{j:k<j\le n\}}}(1-\beta_j)\bigr]\beta_k\rho\nrmB{(\xx_k-\xx^*)}\\
&\le\rho\textstyle\sum_{k=m}^{n}\textstyle\bigl[\prod_{_{\{j:k<j\le n\}}}(1-\beta_j)\bigr]\beta_k\zeta_k \\
&\le \rho\zeta_m \\
\end{split}\end{equation}   
Hence for any $\delta>0$, for all $m$
\begin{equation}\label{R27}
\limsup_{n\rightarrow\infty}\nrmB{\xv_{m,n}}\le 
\rho\zeta_m.\end{equation}
Combining \eqref{R19} along with \eqref{b35}, \eqref{b33} and  \eqref{R27}, we get for any $m\ge m_\delta$ 
\begin{equation}\label{b40}
\limsup_{n\rightarrow\infty}\nrmB{\xx_{n+1}-\xx^*}\le \rho\zeta_m +\delta.\end{equation}
Now taking limit as $m\rightarrow\infty$ on the RHS above and then limit as $\delta\downarrow 0$, we get
\[\limsup_{n\rightarrow\infty}\nrmB{\xx_{n+1}}\le \rho\zeta^*.\]
From the definition of $\zeta^*$, it follows that  $\limsup_{n\rightarrow\infty}\nrmB{\xx_{n+1}-\xx^*}=\zeta^*$ and hence we have shown
\begin{equation}\label{b43}  
\zeta^*\le \rho\zeta^*.
\end{equation}      
Since $\zeta^*\le \gamma^*<\infty$ and $\rho<1$, \eqref{b43} implies that $\zeta^*=0$. It follows that $\limsup_{n\rightarrow\infty}\nrmB{\xx_{n+1}-\xx^*}=0$.\qed
\section{Proofs of Main Results}    
We now come to the proofs of the 4 main theorems. We begin with the first one.\\
{\em Proof of Theorem \ref{T1}}:\\
Let  $\{\CZ_{n+1}:n\ge 0\}$ and  $\{\alpha_n: n\ge 0\}$ be as in Theorem \ref{T1}.
Let $\{\CS_{n+1}:\; n\ge 0\}$ be defined by
\begin{equation}\label{3g}
\CS_{n+1}=\textstyle\sum_{k=0}^n\alpha_k\CZ_{k+1}.
 \end{equation}
We will show that $\{\CS_{n+1}:\; n\ge 0\}$ converges almost surely. Then for any $\omega\in\Omega$ such that $\CS_{n+1}(\omega)$ converges,  using Theorem \ref{T0}, it would follow that $\CX_{n+1}(\omega)$ defined by \eqref{3b} converges to $\xx^*$, completing the proof of Theorem \ref{T1}. 

Note that for $m\le n$, $\CS_{n+1}-\CS_{m}=\sum_{k=m}^n\alpha_k\CZ_{k+1}$. Since  
$\{\CZ_{n+1}:\; n\ge 0\}$ are iid  Gaussian variables, it follows that 
\[\frac{1}{\sqrt{\sum_{k=m}^n\alpha_k^2}}(\CS_{n+1}-\CS_{m})\]
has the same distribution as $\CZ_1$. Thus, Fernique's theorem implies that
\[
\E\bigl[\textstyle\nrmB{\frac{1}{\sqrt{\sum_{k=m}^n\alpha_k^2}}(\CS_{n+1}-\CS_{m})}^2\bigr]=\E\bigl[\nrmB{\CZ_1}^2\bigr]\bigr]=K<\infty.\]
As a consequence,
\[
\E\bigl[\nrmB{(\CS_{n+1}-\CS_{m})}^2\bigr]\le K\textstyle\sum_{k=m}^n\alpha_k^2.\]
In view of assumption \eqref{3a}, it follows  that 
\begin{equation}\label{3h}
\lim_{j\rightarrow\infty}\sup_{ n\ge m\ge j}\E\bigl[\nrmB{(\CS_{n+1}-\CS_{m})}^2\bigr]=0.
\end{equation}
Thus $\CS_{n+1}$ converges in probability (in $\nrmB{\cdot}$). Using the Ito-Nisio theorem we conclude that  $\CS_{n+1}$  converges almost surely completing the proof.\qed

{\em Proof of Theorem \ref{T2}}:\\
Let  $\{\CZ_{n+1}:n\ge 0\}$ and  $\{\alpha_n: n\ge 0\}$ be as in Theorem \ref{T2}.
Let $\{\CS_{n+1}:\; n\ge 0\}$ be defined by \eqref{3g}. We will show that \eqref{3h} holds. The rest of the proof is same as the one for Theorem \ref{T1}.

Let for $0\le t\le 1$, $\clf_t=\sigma(\CZ_{n,s}, \;0\le\le t,\;n\ge 0)$. Since $\CZ_n$ are independent and for each $n$, $\{\CZ_{n,t}: 0\le t\le 1\}$ is a martingale, it follows that
$(\CS_{n,t},\clf_t)_{0\le t\le 1}$ is a martingale for each $n$ and hence so is 
$(\CS_{n+1,t}-\CS_{m,t},\clf_t)_{0\le t\le 1}$. It follows by Doob's maximal inequality that
\[\begin{split}
\E\bigl[\sup_{0\le t\le 1}|\CS_{n+1,t}-\CS_{m,t}|^2\bigr]&\le 4\E\bigl[|\CS_{n+1,1}-\CS_{m,1}|^2\bigr]\\
&\textstyle\le4\sum_{k=m}^n\E\bigl[|\alpha_k\CZ_{k+1,1}|^2\bigr]\\
&\textstyle\le4 \sum_{k=m}^n\alpha_k^2\sigma^2.\end{split}\]
Thus
\begin{equation}\label{3j}
\E\bigl[\nrmB{(\CS_{n+1}-\CS_{m})}^2\bigr]\le \sigma^2\textstyle\sum_{k=m}^n\alpha_k^2.
\end{equation}
Thus \eqref{3h} is true. Since $\{\CZ_k:k\ge 0\}$ are independent,  the conclusion follows using Ito-Nisio theorem as in Theorem \ref{T1}.
\qed

{\em Proof of Theorem \ref{T3}}:\\
Let  $\{\CZ_{n+1}:n\ge 0\}$ and  $\{\alpha_n: n\ge 0\}$ be as in Theorem \ref{T3}. 
Since in Theorem \ref{T0} we required $\{\psi_n: n\ge 0\}$ to be increasing, more precisely satisfy \eqref{R12} while no such condition is put on $\{\lambda_n: n\ge 0\}$ (see \eqref{3b1}),  we proceed as follows. For $n\ge 0$, let
\[\psi_n(\xu_0,\xu_1,\ldots,\xu_n)=C(1+\max_{0\le k\le n}\nrmB{\xu_k})\]
where $C$ is the constant appearing in condition \eqref{3b1} on $\{\lambda_n: n\ge 0\}$ and let
\begin{equation}\label{3b3}\CW_{n+1}=\frac{\lambda_n(\CY_0,\ldots,\CY_n)}{\psi_n(\CY_0,\ldots,\CY_n)}\CZ_{n+1}.\end{equation}
Then \eqref{3b2} can be expressed as
\begin{equation}\label{3b5}
\CY_{n+1}=\CY_n-\alpha_n\Bigl(\Gx(\CY_n)+\psi_n(\CY_0,\ldots,\CY_n)\CW_{n+1}\Bigr),\;\;n\ge 0.\end{equation}
Note that $\{\psi_n:n\ge 0\}$ satisfies \eqref{R12} by definition. Once again, we will show that 
\begin{equation}\label{5a9}
\CS_{n+1}=\textstyle\sum_{k=0}^n\CW_{k+1}\;\;\;\text{converges almost surely in $\nrmB{\cdot}$}.
\end{equation}
{The rest of the proof will be
is as in the proof of Theorem \ref{T1}, once again invoking Theorem \ref{T0}.

For $n\ge 0$ let
\begin{equation}\label{3b4}\begin{split}
\widehat{\CW}_{n+1}&=\CW_{n+1}\Ind_{_{\{\nrmB{\alpha_{n}\CZ_{n+1}}\ge 1\}}}\\
\overline{\CW}_{n+1}&=\E[\CW_{n+1}\Ind_{_{\{\nrmB{\alpha_{n}\CZ_{n+1}}<1\}}}|\sigma(\CZ_j:j\le n)]\\
\widetilde{\CW}_{n+1}&={\CW}_{n+1}\Ind_{\{\nrmB{\alpha_{n}\CZ_{n+1}}<1\}} - \overline{\CW}_{n+1} .\end{split}\end{equation}

Since $0< \chi_n=\frac{\lambda_n}{\psi_n}\le 1$, and $\chi_n$ is $\{\sigma(\CZ_j:j\le n): n\ge 0\}$ adapted, the properties \eqref{4a} and \eqref{4b} of $\{\CZ_n:n\ge 0\}$ yield that for all $n\ge 0$, we also have 
\begin{equation}\label{5a1}
\E\bigl[\nrmB{\overline{\CW}_{n+1}}\bigr]\le \mu_n,
\end{equation}
\begin{equation}\label{5b1}
\E\bigl[\nrmB{\widetilde{\CW}_{n+1}|\sigma(\CZ_j:j\le n)]}^2)\bigr]\le \sigma^2_n.
\end{equation}
Let  $\widehat{\CS}_{n+1}=\sum_{k=0}^n\alpha_k\widehat{\CW}_{k+1}$, 
$\overline{\CS}_{n+1}=\sum_{k=0}^n\alpha_k\overline{\CW}_{k+1}$, 
$\widetilde{\CS}_{n+1}=\sum_{k=0}^n\alpha_k\widetilde{\CW}_{k+1}$.

In view of the assumption \eqref{4}, it follows that 
\begin{equation}\label{5b10}
\textstyle\sum_{k=1}^\infty\P(\widehat{\CW}_{k+1}\neq 0)\le \sum_{k=1}^\infty\P(\widehat{\CZ}_{k+1}\neq 0)<\infty
\end{equation}
and hence using Borel-Cantelli lemma it follows that $\nrmB{\widehat{\CS}_{n+1}}$ converges to zero almost surely.

In view of assumption \eqref{4c1}, it follows that 
\begin{equation}\label{5c1}\begin{split}
\E\bigl[\textstyle \sum_{k=0}^n\alpha_k\nrmB{\overline{\CW}_{k+1}}\bigr]&\le\E\bigl[\textstyle \sum_{k=0}^n\alpha_k\nrmB{\overline{\CZ}_{k+1}}\bigr],\\
&\le\textstyle \sum_{k=0}^\infty\alpha_k\mu_k\\
&<\infty.
\end{split}
\end{equation}
Thus, it follows that $\textstyle \sum_{k=0}^n\alpha_k\nrmB{\overline{\CW}_{k+1}}$ converges almost surely as $n\rightarrow\infty$.
Noting that for $m\le n$
\[
\nrmB{\overline{\CS}_{n+1}-\overline{\CS}_{m+1}}\le\textstyle \sum_{k=m}^n\alpha_k\nrmB{\overline{\CW}_{k+1}}\]
we conclude that $\overline{\CS}_{n+1}$ converges  almost surely in $\nrmB{\cdot}$. 

Since $\widehat{\CS}_{n+1}=\widehat{\CS}_{k+1}+\overline{\CS}_{n+1}+\widetilde{\CS}_{n+1}$,
to complete proof of \eqref{5a9}, remains to prove that $\widetilde{\CS}_{n+1}$ converges in $\nrmB{\cdot}$ almost surely,

By its definition,  $\{\widetilde{\CS}_{n+1}:n\ge 0\}$ is a $\Bb$- valued martingale. Since $\widetilde{\CS}_{n+1}-\widetilde{\CS}_{n}=\alpha_n\widetilde{\CW}_{n+1}$,  \eqref{5b1} implies
\begin{equation}\label{5u}
\E[\nrmB{\widetilde{\CS}_{n+1}-\widetilde{\CS}_{n}}^2)]\le \alpha_n^2\sigma^2_n,\;\;n\ge 0.
\end{equation}
If $1<p<2$, since for any real valued random variable $U$, $(\E[|U|^p])^\frac{p}{2}\le E[|U|^2]$, we conclude that
\begin{equation}\label{5v}
\E[\nrmB{\widetilde{\CS}_{n+1}-\widetilde{\CS}_{n}}^p)]\le \alpha_n^p\sigma^p_n,\;\;n\ge 0
\end{equation}
and hence \eqref{5u} and \eqref{5v} along with the assumptions that $\Bb$ is $p$ uniformly smooth and \eqref{4d1} yields
\begin{equation}\label{5b10}
\textstyle\sum_{n=0}^\infty \E[\nrmB{\widetilde{\CS}_{n+1}-\widetilde{\CS}_{n}}^p)]<\infty.\end{equation}
Theorem 3.2.2 in  Woyczynski \cite{Woyczynski} along with \eqref{5b10} imply that  the martingale $\widetilde{\CS}_{n+1}$ converges in $\nrmB{\cdot}$ almost surely.

This completes the proof.\qed

{\em Proof of Theorem \ref{T4}}:\\
The proof follows steps in the proof of Theorem \ref{T3}.  Defining $\CW_{n+1}$, $\CY_{n+1}$, $\CS_{n+1}$ by \eqref{3b3}, \eqref{3b5}, \eqref{5a9} and we modify \eqref{3b4} as follows:
\begin{equation}\label{3b7}\begin{split}
\widehat{\CW}_{n+1,t}&=\CW_{n+1,t}\Ind_{_{\{\nrm{\alpha_{n}\CZ_{n+1,t}}\ge 1\}}}\\
\overline{\CW}_{n+1,t}&=\E[\CW_{n+1,t}\Ind_{_{\{\nrm{\alpha_{n}\CZ_{n+1,t}}<1\}}}|\sigma(\CW_j:j\le n)]\\
\widetilde{\CW}_{n+1,t}&={\CW}_{n+1,t}\Ind_{_{\{\nrm{\alpha_{n}\CZ_{n+1,t}}<1\}}}-\overline{\CW}_{n+1,t}\Bigr.\end{split}\end{equation}
The equation \eqref{4p} can be rewritten as 
\begin{equation}\label{4px}
\E\bigl[\nrmB{\widehat{\CW}_{n+1}}\bigr]\le\delta_n,
\end{equation}
and the equation \eqref{4q} can be written as 
\begin{equation}\label{4qx}
\E\bigl[\nrmB{\overline{\CW}_{n+1}}]\le\mu_n,
\end{equation}
Defining   $\widehat{\CS}_{n+1,t}=\sum_{k=0}^n\alpha_k\widehat{\CW}_{k+1,t}$, 
$\overline{\CS}_{n+1,t}=\sum_{k=0}^n\alpha_k\overline{\CW}_{k+1,t}$, 
$\widetilde{\CS}_{n+1,t}=\sum_{k=0}^n\alpha_k\widetilde{\CW}_{k+1,t}$ and using  \eqref{4px} and \eqref{4qx} we now have
\begin{equation}\label{5b12}
\E\textstyle\bigl[\sum_{k=1}^\infty\alpha_k\nrmB{\widehat{\CW}_{k+1}}\bigr]<\infty ,\end{equation}
\begin{equation}\label{5b13}
\E\textstyle\bigl[\sum_{k=1}^\infty\alpha_k\nrmB{\overline{\CW}_{k+1}}\bigr]<\infty.\end{equation}
 Since for $m< n$
\[\textstyle
\bigl[\nrmB{\widehat{\CS}_{n+1}-\widehat{\CS}_{m+1}}\bigr]\le \sum_{k=m+1}^n\alpha_k\bigl[\nrmB{\widehat{\CW}_{k+1}}\bigr],
\]
it follows from \eqref{5b12} that $\widehat{\CS}_{n+1}$ converges almost surely. Similarly, using  \eqref{5b13} we can conclude that $\overline{\CS}_{n+1}$ converges almost surely. 
The equation \eqref{4r} gives us
\begin{equation}\label{4rx}
\textstyle\E\bigl[\int_0^1\bigl|\widetilde{\CW}_{k+1,t}|^2dt\bigr]\le \sigma^2_n,\;\;n\ge 0.
\end{equation}
Since $|\widetilde{\CW}_{k+1,t}|$ is bounded, it follows that $\{\widetilde{\CW}_{k+1}:k\ge 0\}$ is $\L^q([0,1],\R^d)$ valued for every $q\in[1,\infty)$. 
By its definition, it is clear that $\{\widetilde{\CS}_{k+1}:k\ge 0\}$ is a martingale and \eqref{4rx} yields
\begin{equation}\label{4ry}
\begin{split}
\textstyle\E\bigl[\int_0^1\bigl|\widetilde{\CS}_{n+1,t}|^2dt\bigr]
&\le \E\bigl[\textstyle\sum_{k=1}^n\alpha_k^2\int_0^1\bigl|\widetilde{\CW}_{k+1,t}|^2dt\bigr]\\
&\le \textstyle\sum_{k=1}^n\alpha_k^2\sigma^2_k\\
&\le \textstyle\sum_{k=1}^\infty\alpha_k^2\sigma^2_k\\
&<\infty.\end{split}
\end{equation}
Thus $\{\widetilde{\CS}_{n+1}:n\ge 0\}$ is an $\L^2([0,1],\R^d)$ valued martingale with bounded second moments. Thus $\widetilde{\CS}_{n+1}$ converges almost surely as $n\rightarrow\infty$ in $\L^2([0,1],\R^d)$. If $p<2$, then $\L^2([0,1],\R^d)$ convergence implies $\L^p([0,1],\R^d)$. This completes the proof of Theorem \ref{T4}.  \qed
\begin{remark}
The $\L^p([0,1],\R^d)$ can be replaced by $\L^p((\Gamma,\cla,\pi),\R^d)$, where $\pi$ is a probability measure on a measurable space $(\Gamma,\cla)$ - same proof works.
\end{remark}

\end{document}